\documentclass[a4paper,11pt]{amsart}

\usepackage{graphicx}
\usepackage{mathptmx}
\usepackage{amsmath}
\usepackage{amssymb}
\usepackage{enumitem}
\usepackage{xcolor}

\newmuskip\pFqmuskip

\newcommand*\pFq[6][8]{%
  \begingroup 
  \pFqmuskip=#1mu\relax
  \mathcode`=\string"8000
  \begingroup\lccode`\~=`\,
  \lowercase{\endgroup\let~}\pFqcomma
  F^{#2}_{#3}{\left(\genfrac..{0pt}{}{#4}{#5}\bigg|#6\right)}%
  \endgroup
}
\newcommand{\pFqcomma}{\mskip\pFqmuskip}

\newtheorem{theorem}{Theorem}[section]

\begin{document}

\title[Probabilistic bivariate Bell polynomials]{Probabilistic bivariate Bell polynomials}

\author{Taekyun  Kim}
\address{Department of Mathematics, Kwangwoon University, Seoul 139-701, Republic of Korea}
\email{tkkim@kw.ac.kr}
\author{Dae San  Kim}
\address{Department of Mathematics, Sogang University, Seoul 121-742, Republic of Korea}
\email{dskim@sogang.ac.kr}


\maketitle

\begin{abstract}
Let $Y$ be a random variable whose moment generating function exists in some neighborhood of the origin. We consider the probabilistic bivariate Bell polynomials associated with $Y$ and the probabilistic bivariate $r$-Bell polynomials associated with $Y$. For those polynomials, we derive the recurrence relations corresponding to the ones found by Zheng and Li for the bivariate Bell polynomials and the bivariate $r$-Bell polynomials.
\end{abstract}

\section{Introduction}
The Stirling number of the second kind ${n\brace k}$ enumerates the number of ways to partition a set of $n$ objects into $k$ nonempty subsets (see \eqref{1}). More generally, for any $r\in\mathbb{N}$, the $r$-Stirling number of the second kind ${n+r \brace k+r}_{r}$ counts the number of ways to partition a set of $n+r$ objects into $k+r$ nonempty subsets such that the first $r$ elements are in distinct subsets (see \eqref{2}). \par
Spivey discovered a remarkable recurrence relation for the Bell numbers (see \eqref{4}). This relation was extended to the Bell polynomials by Gould and Quaintance (see \eqref{5}). The corresponding ones for the bivariate Bell polynomials $\phi_{n}(x,y)=\sum_{k=0}^{n}{n\brace k}(x)_{k}y^{k}$ and the bivariate $r$-Bell polynomials $\phi_{n,r}(x,y)=\sum_{k=0}^{n}{n+r\brace k+r}_{r}(x)_{k}y^{k}$
were found by Zheng and Li (see \eqref{9}, \eqref{10}).\par
Assume that $Y$ is a random variable whose moment generating function exists in some neighbourhood of the origin. Let $(Y_{j})_{j\ge 1}$ be a sequence of mutually independent copies of the random variable $Y$, and let $S_{k}=Y_{1}+Y_{2}+\cdots+Y_{k},\ (k\ge 1)$, with $S_{0}=0$. In this paper, we introduce bivariate extensions of the probabilistic Bell polynomials associated $Y$ and the probabilistic $r$-Bell polynomials associated with $Y$, namely the probabilistic bivariate Bell polynomials associated with $Y$, $\phi_{n}^{Y}(x,y)$ (see \eqref{16}) and the probabilistic bivariate $r$-Bell polynomials associated with $Y$, $\phi_{n,r}^{Y}(x,y)$ (see \eqref{20}). Then we derive the recurrence relations corresponding to the ones found by Zheng and Li for the probabilistic bivariate Bell polynomials associated with $Y$ (see Theorem 2.2)
\begin{align*}
&\phi_{m+n}^{Y}(x,y)\\
&=\sum_{k=0}^{n}\sum_{j=0}^{m}\binom{m}{j}\phi_{j}^{Y}(x-k,y)(x)_{k}y^{k}\frac{1}{k!}\sum_{l_{1}+\cdots+l_{k}=n}\binom{n}{l_{1},\dots,l_{k}}E\bigg[S_{k}^{m-j}\prod_{i=1}^{k}Y_{i}^{l_{i}}\bigg],
\end{align*}
and also for the probabilistic bivariate $r$-Bell polynomials associated with $Y$ (see Theorems 2.4 and 2.5)
\begin{align*}
&\phi_{m+n,r}^{Y}(x,y)\\
&= \sum_{i=0}^{n}\sum_{k=0}^{m}\binom{n}{i} \phi_{i,r}^{Y}(x-k,y) \frac{(x)_{k}y^{k}}{k!} \sum_{j=k}^{m}\sum_{l_{1}+\cdots+l_{k}=j}\binom{m}{j}\binom{j}{l_{1},\dots,l_{k}}r^{m-j} E\bigg[S_{k}^{n-i}\prod_{i=1}^{k}Y_{i}^{l_{i}}\bigg] \\
&= \sum_{i=0}^{n}\sum_{k=0}^{m}\frac{1}{k!}\binom{n}{i}(x)_{k}y^{k}\phi_{i}^{Y}(x-k,y) \sum_{j=k}^{m}\binom{m}{j}\sum_{l_{1}+\cdots+l_{k}=j}\binom{j}{l_{1},\dots,l_{k}} E\bigg[(S_{k}+r)^{n-i}\prod_{i=1}^{k}Y_{i}^{l_{i}}\bigg] r^{m-j}.
\end{align*}
These are done by using the simple formula for Taylor series (see \eqref{11}). For the rest of this section, we recall the facts that are needed throughout this paper. \par

\vspace{0.1in}

For $n\ge 0$, the Stirling numbers of the second kind are given by
\begin{equation}
x^{n}=\sum_{k=0}^{n}{n\brace k}(x)_{k},\quad (\mathrm{see}\ [1-29]),\label{1}	
\end{equation}
where $(x)_{0}=1,\ (x)_{k}=x(x-1)\cdots(x-k+1),\ (k\ge 1)$. \\
The $r$-Stirling number of the second kind are defined as
\begin{equation}
(x+r)^{n}=\sum_{k=0}^{n}{n+r \brace k+r}_{r}(x)_{k},\quad (n\ge 0),\quad (\mathrm{see}\ [8,18-20]).\label{2}
\end{equation}
The Bell polynomials are given by
\begin{equation}
\phi_{n}(x)=\sum_{k=0}^{n}{n\brace k}x^{k},\quad (n\ge 0),\quad (\mathrm{see}\ [2,10,17,25,28]).\label{3}
\end{equation}
When $x=1,\ \phi_{n}=\phi_{n}(1)$ are called the Bell numbers. \par
Spivey found the following recurrence relation for $\phi_{n}$:
\begin{equation}
\phi_{l+n}=\sum_{k=0}^{l}\sum_{i=0}^{n}k^{n-i}\binom{n}{i}{l \brace k}\phi_{i},\quad (\mathrm{see}\ [28]).\label{4}
\end{equation} \\
Gould and Quaintance discovered the following extension of the relation for $\phi_{n}$ in \eqref{4} to that for $\phi_{n}(x)$:
\begin{equation}
\phi_{l+n}(x)=\sum_{k=0}^{l}\sum_{i=0}^{n}k^{n-i}\binom{n}{i}{l \brace k}\phi_{i}(x)x^{k},\quad (\mathrm{see}\ [12,28]).\label{5}	
\end{equation} \par
For $r\in\mathbb{N}$, the $r$-Bell polynomials are defined as
\begin{equation}
\phi_{n,r}(x)=\sum_{k=0}^{n}{n+r \brace k+r}_{r}x^{k},\quad (n\ge 0),\quad (\mathrm{see}\ [18,20]).\label{6}
\end{equation}
It is known that the bivariate Bell polynomials are given by
\begin{equation}
\phi_{n}(x,y)=\sum_{k=0}^{n}{n\brace k}(x)_{k}y^{k},\quad (\mathrm{see}\ [30]).\label{7}	
\end{equation}
For $r \in \mathbb{N}$, the bivariate $r$-Bell polynomials are defined as
\begin{equation}
\phi_{n,r}(x,y)=\sum_{k=0}^{n}{n+r\brace k+r}_{r}(x)_{k}y^{k},\quad (\mathrm{see}\ [30]).\label{8}
\end{equation}
Note that
\begin{displaymath}
\lim_{y\rightarrow 0}\phi_{n}\bigg(\frac{x}{y},y\bigg)=\phi_{n}(x),\quad
\lim_{y\rightarrow 0}\phi_{n,r}\bigg(\frac{x}{y},y\bigg)=\phi_{n,r}(x),\quad (n\ge 0).
\end{displaymath}
From \eqref{7} and \eqref{8}, we note that
\begin{align}
\phi_{m+n}(x,y)&=\sum_{k=0}^{m}\sum_{i=0}^{n}k^{n-i}\binom{n}{i}{m\brace k}\phi_{i}(x-k,y)(x)_{k}y^{k}, \label{9} 	\\
\phi_{m+n,r}(x,y)&=\sum_{k=0}^{m}\sum_{i=0}^{n}k^{n-i}\binom{n}{i}{m+r \brace k+r}_{r}\phi_{i,r}(x-k,y)(x)_{k}y^{k},\quad (\mathrm{see}\ [30]), \label{10}
\end{align}
where $m,n$ are nonnegative integers. \par
We observe that
\begin{align}
f(x+t)&=\sum_{n=0}^{\infty}\frac{f^{(n)}(x)}{n!}(x+t-x)^{n}=\sum_{n=0}^{\infty}\frac{f^{(n)}(x)}{n!}t^{n}\label{11} \\
&=\sum_{n=0}^{\infty}\frac{t^{n}D_{x}^{n}}{n!}f(x)=e^{tD_{x}}f(x), \nonumber
\end{align}
where $D_{x}=\frac{d}{dx}$. \par
Assume that $Y$ is a random variable such that the moment generating function of $Y$,
\begin{equation*}
E\big[e^{tY}\big]=\sum_{n=0}^{\infty}\frac{t^{n}}{n!}E[Y^{n}],\quad (|t|<r),
\end{equation*}
exists for some $r>0$, (see [3,16,17]).
Let $(Y_{j})_{j\ge 1}$ be a sequence of mutually independent copies of the random variable $Y$, and let $S_{k}=Y_{1}+Y_{2}+\cdots+Y_{k},\ (k\ge 1)$, with $S_{0}=0$. \par
We recall that the probabilistic Stirling numbers of the second kind associated with $Y$ are given by
\begin{equation}
{n\brace k}_{Y}=\frac{1}{k!}\sum_{i=0}^{k}\binom{k}{i}(-1)^{k-i}E\big[S_{i}^{n}\big],\quad (0\le k\le n),\quad (\mathrm{see}\ [3,16,17]). \label{12}
\end{equation}
The probabilistic Bell polynomials associated with $Y$ are defined as
\begin{equation}
\phi_{n}^{Y}(x)=\sum_{k=0}^{n}{n\brace k}_{Y}x^{k},\quad (n\ge 0),\quad (\mathrm{see}\ [16,17,27]).\label{13}
\end{equation}
For $r\in\mathbb{N}$, the probabilistic $r$-Stirling numbers of the second kind associated with $Y$ are given by
\begin{equation}
{n+r \brace k+r}_{r,Y}=\frac{1}{k!}\sum_{j=0}^{k}\binom{k}{j}(-1)^{k-j}E\big[(S_{j}+r)^{n}\big],\quad (n\ge k\ge 0),\quad (\mathrm{see}\ [16,17]). \label{14}
\end{equation}
When $Y=1$, we have
\begin{displaymath}
{n+r\brace k+r}_{r,Y}={n+r \brace k+r}_{r},\ (n\ge k\ge 0).
\end{displaymath}
The probabilistic $r$-Bell polynomials associated with $Y$ are defined as
\begin{equation}
\phi_{n,r}^{Y}(x)=\sum_{k=0}^{n}{n+r \brace k+r}_{r,Y}x^{k},\quad (n\ge 0),\quad (\mathrm{see}\ [16,17]).\label{15}
\end{equation}
When $Y=1$, we have $\phi_{n,r}^{Y}(x)=\phi_{n,r}(x),\ (n\ge 0)$. \par

\section{Probabilistic bivariate Bell polynomials}
Let $(Y_{j})_{j\ge 1}$ be a sequence of mutually independent copies of the random variable $Y$, and let
\begin{displaymath}
S_{0}=0,\quad S_{k}=Y_{1}+Y_{2}+\cdots+Y_{k},\ (k\ge 1).
\end{displaymath} \par
Now, we define the {\it{probabilistic bivariate Bell polynomials associated with $Y$}} by
\begin{equation}
\Big(1+y\big(E\big[e^{Yt}\big]-1\big)\Big)^{x}=\sum_{n=0}^{\infty}\phi_{n}^{Y}(x,y)\frac{t^{n}}{n!}.\label{16}
\end{equation}
From \eqref{16}, we note that
\begin{align}
\sum_{n=0}^{\infty}&\phi_{n}^{Y}(x,y)\frac{t^{n}}{n!}=\sum_{k=0}^{\infty}(x)_{k}y^{k}\frac{1}{k!}\Big(E\big[e^{Yt}\big]-1\Big)^{k}\label{17} \\
&=\sum_{k=0}^{\infty}(x)_{k}y^{k}\sum_{n=k}^{\infty}{n \brace k}_{Y}\frac{t^{n}}{n!}=\sum_{n=0}^{\infty}\bigg(\sum_{k=0}^{n}{n \brace k}_{Y}(x)_{k}y^{k}\bigg)\frac{t^{n}}{n!}.\nonumber	
\end{align}
Therefore, by comparing the coefficients on both sides of \eqref{17}, we obtain the following theorem.
\begin{theorem}
For $n\ge 0$, we have
\begin{displaymath}
\phi_{n}^{Y}(x,y)= \sum_{k=0}^{n}{n \brace k}_{Y}(x)_{k}y^{k}.
\end{displaymath}
\end{theorem}
When $Y=1, \phi_{n}^{Y}(x,y)=\phi_{n}(x,y)$.
As $\lim_{y\rightarrow 0}\phi_{n}^{Y}\big(\frac{x}{y},y\big)=\phi_{n}^{Y}(x)$, the probabilistic bivariate Bell polynomials $\phi_{n}^{Y}(x,y)$ are also a bivariate extension of the probabilistic Bell polynomials $\phi_{n}^{Y}(x)$.

From \eqref{16}, we note that
\begin{align}
e^{sD_{t}} \Big(1+y\big(E\big[e^{Yt}\big]-1\big)\Big)^{x}&=\sum_{n=0}^{\infty}\frac{s^{n}}{n!}D_{t}^{n}\sum_{m=0}^{\infty}\phi_{m}^{Y}(x,y)\frac{t^{m}}{m!}\label{18} \\
&=\sum_{n=0}^{\infty}\sum_{m=0}^{\infty}\phi_{m+n}^{Y}(x,y)\frac{s^{n}}{n!}\frac{t^{m}}{m!}.\nonumber
\end{align} \par
On the other hand, by \eqref{11}, we get
\begin{align}
&e^{sD_{t}} \Big(1+y\big(E\big[e^{Yt}\big]-1\big)\Big)^{x}= \Big(1+y\big(E\big[e^{(s+t)Y}\big]-1\big)\Big)^{x}\label{19}\\
&= \Big(1+y\big(E\big[e^{Yt}\big]-1\big)+ yE\big[e^{Yt}(e^{sY}-1)\big]\Big)^{x}\nonumber\\
&=\sum_{k=0}^{\infty}\binom{x}{k} \Big(1+y\big(E\big[e^{Yt}\big]-1\big)\Big)^{x-k}y^{k}\Big(E\big[e^{Yt}(e^{sY}-1)\big]\Big)^{k} \nonumber \\
&=\sum_{k=0}^{\infty}\frac{(x)_{k}}{k!}y^{k}\sum_{j=0}^{\infty}\phi_{j}^{Y}(x-k,y)\frac{t^{j}}{j!}E\Big[e^{S_{k}t}(e^{sY_{1}}-1) (e^{sY_{2}}-1)\cdots (e^{sY_{k}}-1)\Big]\nonumber \\
&=\sum_{k=0}^{\infty}\frac{(x)_{k}}{k!}y^{k}\sum_{n=k}^{\infty}\sum_{l_{1}+l_{2}+\cdots+l_{k}=n}\binom{n}{l_{1},l_{2},\dots,l_{k}}E\Big[e^{S_{k}t}Y_{1}^{l_{1}} Y_{2}^{l_{2}}\cdots Y_{k}^{l_{k}}\Big]\frac{s^{n}}{n!}\sum_{j=0}^{\infty}\phi_{j}^{Y}(x-k,y)\frac{t^{j}}{j!}\nonumber \\
&=\sum_{n=0}^{\infty}\sum_{k=0}^{n}(x)_{k}y^{k}\frac{1}{k!}\sum_{l_{1}+\cdots+l_{k}=n}\binom{n}{l_{1},\dots,l_{k}}\frac{s^{n}}{n!}\sum_{l=0}^{\infty}E\bigg[S_{k}^{l}\prod_{i=1}^{k}Y_{i}^{l_{i}}\bigg]\frac{t^{l}}{l!}\sum_{j=0}^{\infty}\phi_{j}^{Y}(x-k,y)\frac{t^{j}}{j!}\nonumber \\
&=\sum_{n=0}^{\infty}\sum_{m=0}^{\infty}\sum_{k=0}^{n}\sum_{j=0}^{m}\binom{m}{j}\phi_{j}^{Y}(x-k,y)(x)_{k}y^{k}\frac{1}{k!}\sum_{l_{1}+\cdots+l_{k}=n}\binom{n}{l_{1},\dots,l_{k}}E\bigg[S_{k}^{m-j}\prod_{i=1}^{k}Y_{i}^{l_{i}}\bigg]\frac{s^{n}}{n!}\frac{t^{m}}{m!}, \nonumber
\end{align}
where $l_{1},l_{2},\dots,l_{k}$ are positive integers. \\
Therefore, by \eqref{18} and \eqref{19}, we obtain the following theorem.
\begin{theorem}
For $m,n\ge 0$, we have
\begin{align*}
&\phi_{m+n}^{Y}(x,y)\\
&=\sum_{k=0}^{n}\sum_{j=0}^{m}\binom{m}{j}\phi_{j}^{Y}(x-k,y)(x)_{k}y^{k}\frac{1}{k!}\sum_{l_{1}+\cdots+l_{k}=n}\binom{n}{l_{1},\dots,l_{k}}E\bigg[S_{k}^{m-j}\prod_{i=1}^{k}Y_{i}^{l_{i}}\bigg],
\end{align*}
where $l_{1},l_{2},\dots,l_{k}$ are positive integers.
\end{theorem}
When $Y=1$, the expression in Theorem 2.2 reduces to
\begin{align*}
\phi_{m+n}(x,y)&= \sum_{k=0}^{n}\sum_{j=0}^{m}\binom{m}{j}\phi_{j}(x-k,y)(x)_{k}y^{k}\frac{1}{k!}\sum_{l_{1}+\cdots+l_{k}=n}\binom{n}{l_{1},\dots,l_{k}}k^{m-j} \nonumber \\
&= \sum_{k=0}^{n}\sum_{j=0}^{m}\binom{m}{j}\phi_{j}(x-k,y)(x)_{k}y^{k}{n\brace k}k^{m-j}.
\end{align*} \par
For $r\in\mathbb{N}$, we define the {\it{probabilistic bivariate $r$-Bell polynomials associated with $Y$}} by
\begin{equation}
\Big(1+y\big(E\big[e^{Yt}\big]-1\big)\Big)^{x}e^{rt}=\sum_{n=0}^{\infty}\phi_{n,r}^{Y}(x,y)\frac{t^{n}}{n!}.\label{20}	
\end{equation}
From \eqref{14}, we note that
\begin{align}
\sum_{n=0}^{\infty}\phi_{n,r}^{Y}(x,y)\frac{t^{n}}{n!}&= \Big(1+y\big(E\big[e^{Yt}\big]-1\big)\Big)^{x}e^{rt}=\sum_{k=0}^{\infty}(x)_{k}\frac{y^{k}}{k!} \big(E\big[e^{Yt}\big]-1\big)^{k}e^{rt}\label{21} \\
&=\sum_{k=0}^{\infty}y^{k}(x)_{k}\sum_{n=k}^{\infty}{n+r \brace k+r}_{r,Y}\frac{t^{n}}{n!}=\sum_{n=0}^{\infty}\sum_{k=0}^{n}y^{k}(x)_{k}{n+r \brace k+r}_{r,Y}\frac{t^{n}}{n!}.\nonumber
\end{align}
Therefore, by comparing the coefficients on both sides of \eqref{21}, we obtain the following theorem.
\begin{theorem}
For $n\ge 0$, we have
\begin{displaymath}
\phi_{n,r}^{Y}(x,y)= \sum_{k=0}^{n}{n+r \brace k+r}_{r,Y}(x)_{k}y^{k}.
\end{displaymath}
\end{theorem}
When $Y=1$, we obtain
\begin{displaymath}
\phi_{n,r}(x,y)=\sum_{k=0}^{n}{n+r \brace k+r}_{r}(x)_{k}y^{k}.
\end{displaymath}
In addition, we note that
\begin{displaymath}
\lim_{y\rightarrow 0}\phi_{n,r}^{Y}\bigg(\frac{x}{y},y\bigg)=\sum_{k=0}^{n}{n+r \brace k+r}_{r,Y}x^{k}=\phi_{n,r}^{Y}(x).
\end{displaymath} \par
By using \eqref{11}, we have
\begin{align}
&e^{tD_{s}} \Big(1+y\big(E\big[e^{sY}\big]-1\big)\Big)^{x}	e^{rs}= \Big(1+y\big(E\big[e^{(s+t)Y}\big]-1\big)\Big)^{x}e^{(s+t)r}\label{22} \\
&\sum_{k=0}^{\infty}\binom{x}{k} \Big(1+y\big(E\big[e^{tY}\big]-1\big)\Big)^{x-k}e^{rt}y^{k}\Big(E\big[e^{tY}\big(e^{sY}-1\big)\big]\Big)^{k}e^{sr}\nonumber \\
&=\sum_{k=0}^{\infty}\frac{(x)_{k}y^{k}}{k!}\sum_{i=0}^{\infty}\phi_{i,r}^{Y}(x-k,y)\frac{t^{i}}{i!}\sum_{j=k}^{\infty}\sum_{l_{1}+\cdots+l_{k}=j}\binom{j}{l_{1},\dots,l_{k}}E\bigg[e^{S_{k}t}\prod_{i=1}^{k}Y_{i}^{l_{i}}\bigg]\frac{s^{j}}{j!}e^{sr} \nonumber \\
&=\sum_{k=0}^{\infty}\frac{(x)_{k}y^{k}}{k!}\sum_{m=k}^{\infty}\sum_{j=k}^{m}\binom{m}{j} \sum_{l_{1}+\cdots+l_{k}=j}\binom{j}{l_{1},\dots,l_{k}}E\bigg[e^{S_{k}t}\prod_{i=1}^{k}Y_{i}^{l_{i}}\bigg]r^{m-j}\frac{s^{m}}{m!}\sum_{i=0}^{\infty}\phi_{i,r}^{Y}(x-k,y)\frac{t^{i}}{i!}\nonumber \\
&=\sum_{m=0}^{\infty}\sum_{l=0}^{\infty}\sum_{k=0}^{m}\frac{(x)_{k}y^{k}}{k!}\sum_{j=k}^{m}\sum_{l_{1}+\cdots+l_{k}=j}\binom{m}{j}\binom{j}{l_{1},\dots,l_{k}}r^{m-j}\frac{s^{m}}{m!} E\bigg[S_{k}^{l}\prod_{i=1}^{k}Y_{i}^{l_{i}}\bigg]\frac{t^{l}}{l!}\sum_{i=0}^{\infty}\phi_{i,r}^{Y}(x-k,y)\frac{t^{i}}{i!}\nonumber \\
&=\sum_{m=0}^{\infty}\sum_{n=0}^{\infty}\sum_{i=0}^{n}\sum_{k=0}^{m}\binom{n}{i} \phi_{i,r}^{Y}(x-k,y) \frac{(x)_{k}y^{k}}{k!}\nonumber\\
& \quad\quad\quad\quad \times \sum_{j=k}^{m}\sum_{l_{1}+\cdots+l_{k}=j}\binom{m}{j}\binom{j}{l_{1},\dots,l_{k}}r^{m-j} E\bigg[S_{k}^{n-i}\prod_{i=1}^{k}Y_{i}^{l_{i}}\bigg]\frac{t^{n}}{n!}\frac{s^{m}}{m!},\nonumber
\end{align}
where $l_{1},l_{2},\dots,l_{k}$ are positive integers. \par
On the other hand, by \eqref{20}, we get
 \begin{align}
 e^{tD_{s}} \Big(1+y\big(E\big[e^{sY}\big]-1\big)\Big)^{x}	e^{rs}&=\sum_{n=0}^{\infty}\frac{t^{n}}{n!}D_{s}^{n}\sum_{m=0}^{\infty}\phi_{m,r}^{Y}(x,y)\frac{s^{m}}{m!}\label{23} \\
 &=\sum_{n=0}^{\infty}\sum_{m=0}^{\infty}\phi_{m+n,r}^{Y}(x,y)\frac{t^{n}}{n!}\frac{s^{m}}{m!}.\nonumber
 \end{align}
Therefore, by \eqref{22} and \eqref{23}, we obtain the following theorem.
\begin{theorem}
For $m,n\ge 0$, we have
\begin{align*}
&\phi_{m+n,r}^{Y}(x,y)\\
&= \sum_{i=0}^{n}\sum_{k=0}^{m}\binom{n}{i} \phi_{i,r}^{Y}(x-k,y) \frac{(x)_{k}y^{k}}{k!} \sum_{j=k}^{m}\sum_{l_{1}+\cdots+l_{k}=j}\binom{m}{j}\binom{j}{l_{1},\dots,l_{k}}r^{m-j} E\bigg[S_{k}^{n-i}\prod_{i=1}^{k}Y_{i}^{l_{i}}\bigg],
\end{align*}
where $l_{1},l_{2},\dots,l_{k}$ are positive integers.
\end{theorem}
When $Y=1$, the expression in Theorem 2.4 boils down to
\begin{align*}
\phi_{m+n,r}(x,y)&= \sum_{i=0}^{n}\sum_{k=0}^{m}\binom{n}{i} \phi_{i,r}(x-k,y) \frac{(x)_{k}y^{k}}{k!} \sum_{j=k}^{m} \binom{m}{j}\sum_{l_{1}+\cdots+l_{k}=j}\binom{j}{l_{1},\dots,l_{k}}r^{m-j}k^{n-i} \\
&= \sum_{i=0}^{n}\sum_{k=0}^{m}\binom{n}{i} \phi_{i,r}(x-k,y)(x)_{k}y^{k}k^{n-i}{m+r \brace k+r}_{r}.
\end{align*} \par
Now, we observe that
\begin{align}
&e^{tD_{s}} \Big(1+y\big(E\big[e^{sY}\big]-1\big)\Big)^{x}e^{rs}= \Big(1+y\big(E\big[e^{(s+t)Y}\big]-1\big)\Big)^{x}e^{r(s+t)} \label{24} \\
 &= \sum_{k=0}^{\infty}\frac{(x)_{k}y^{k}}{k!} \Big(1+y\big(E\big[e^{tY}\big]-1\big)\Big)^{x-k} E\bigg[e^{(S_{k}+r)t}\prod_{i=1}^{k}\big(e^{sY_{i}}-1\big)\bigg]e^{sr}\nonumber \\
 &= \sum_{k=0}^{\infty}\frac{(x)_{k}y^{k}}{k!}\sum_{i=0}^{\infty}\phi_{i}^{Y}(x-k,y)\frac{t^{i}}{i!}\sum_{j=k}^{\infty}\sum_{l_{1}+\cdots+l_{k}=j}\binom{j}{l_{1},\dots,l_{k}}E\bigg[e^{(S_{k}+r)t}\prod_{i=1}^{k}Y_{i}^{l_{i}}\bigg]\frac{s^{j}}{j!}e^{sr} \nonumber \\
 &=\sum_{k=0}^{\infty}\frac{(x)_{k}y^{k}}{k!}\sum_{m=k}^{\infty}\sum_{j=k}^{m}\binom{m}{j}\sum_{l_{1}+\cdots+l_{k}=j}\binom{j}{l_{1},\dots,l_{k}} E\bigg[e^{(S_{k}+r)t}\prod_{i=1}^{k}Y_{i}^{l_{i}}\bigg] r^{m-j}\frac{s^{m}}{m!}\sum_{i=0}^{\infty}\phi_{i}^{Y}(x-k,y)\frac{t^{i}}{i!} \nonumber \\ 	
 &=\sum_{m=0}^{\infty}\sum_{k=0}^{m} \frac{(x)_{k}y^{k}}{k!}\sum_{j=k}^{m}\binom{m}{j}\sum_{l_{1}+\cdots+l_{k}=j}\binom{j}{l_{1},\dots,l_{k}} E\bigg[e^{(S_{k}+r)t}\prod_{i=1}^{k}Y_{i}^{l_{i}}\bigg] r^{m-j}\frac{s^{m}}{m!}\sum_{i=0}^{\infty}\phi_{i}^{Y}(x-k,y)\frac{t^{i}}{i!}\nonumber \\
 &=\sum_{m=0}^{\infty}\sum_{l=0}^{\infty} \sum_{k=0}^{m}\frac{(x)_{k}y^{k}}{k!}\sum_{j=k}^{m}\binom{m}{j}\sum_{l_{1}+\cdots+l_{k}=j}\binom{j}{l_{1},\dots,l_{k}} r^{m-j}\frac{s^{m}}{m!}\nonumber\\
&\quad\quad\quad\quad \times E\bigg[(S_{k}+r)^{l}\prod_{i=1}^{k}Y_{i}^{l_{i}}\bigg] \frac{t^{l}}{l!} \sum_{i=0}^{\infty}\phi_{i}^{Y}(x-k,y)\frac{t^{i}}{i!}\nonumber\\
 &=\sum_{m=0}^{\infty}\sum_{n=0}^{\infty}\sum_{i=0}^{n}\sum_{k=0}^{m}\binom{n}{i}\frac{(x)_{k}y^{k}}{k!}\phi_{i}^{Y}(x-k,y) \sum_{j=k}^{m}\binom{m}{j}\sum_{l_{1}+\cdots+l_{k}=j}\binom{j}{l_{1},\dots,l_{k}}\nonumber\\
&\quad\quad\quad\quad \times E\bigg[(S_{k}+r)^{n-i}\prod_{i=1}^{k}Y_{i}^{l_{i}}\bigg] r^{m-j}\frac{s^{m}}{m!}\frac{t^{n}}{n!},\nonumber
\end{align}
where $l_{1},l_{2},\dots,l_{k}$ are positive integers. \\
Therefore, by \eqref{23} and \eqref{24}, we obtain the following theorem.
\begin{theorem}
For $m,n\ge 0$, we have
\begin{align*}
&\phi_{m+n,r}^{Y}(x,y)\\
&= \sum_{i=0}^{n}\sum_{k=0}^{m}\frac{1}{k!}\binom{n}{i}(x)_{k}y^{k}\phi_{i}^{Y}(x-k,y) \sum_{j=k}^{m}\binom{m}{j}\sum_{l_{1}+\cdots+l_{k}=j}\binom{j}{l_{1},\dots,l_{k}} E\bigg[(S_{k}+r)^{n-i}\prod_{i=1}^{k}Y_{i}^{l_{i}}\bigg] r^{m-j},
\end{align*}
where $l_{1},l_{2},\dots,l_{k}$ are positive integers.
\end{theorem}
When $Y=1$, the expression in Theorem 2.5 reduces to
\begin{align*}
&\phi_{m+n,r}(x,y)\\
&=\sum_{i=0}^{n}\sum_{k=0}^{m}\binom{n}{i}(x)_{k}y^{k}\phi_{i}(x-k,y)\frac{1}{k!}\sum_{j=k}^{m}\binom{m}{j}\sum_{l_{1}+\cdots+l_{k}=j}\binom{j}{l_{1},\dots,l_{k}}r^{m-j}(r+k)^{n-i}\\
&=\sum_{i=0}^{n}\sum_{k=0}^{m}\binom{n}{i}(x)_{k}y^{k}\phi_{i}(x-k,y){m+r \brace k+r}_{r}(r+k)^{n-i}.
\end{align*} \par
From Theorem 2.4, we note that
\begin{align}
&\phi_{m+n,r}^{Y}(x)=\lim_{y\rightarrow 0}\phi_{m+n}^{Y}\bigg(\frac{x}{y},y\bigg)\label{25} \\
&=\lim_{y\rightarrow 0}\sum_{i=0}^{n}\sum_{k=0}^{m}\binom{n}{i}\phi_{i,r}^{Y}\bigg(\frac{x}{y}-k,y\bigg)\bigg(\frac{x}{y}\bigg)_{k}y^{k}\frac{1}{k!}\nonumber\\
&\quad\quad\quad\quad \times \sum_{j=k}^{m}\binom{m}{j}\sum_{l_{1}+\cdots+l_{k}=j}\binom{j}{l_{1},\dots,l_{k}}r^{m-j}E\bigg[S_{k}^{n-i}\prod_{i=1}^{k}Y_{i}^{l_{i}}\bigg] \nonumber \\
&=\sum_{i=0}^{n}\sum_{k=0}^{m}\binom{n}{i}\phi_{i,r}^{Y}(x)x^{k}\frac{1}{k!}\sum_{j=k}^{m}\binom{m}{j} \sum_{l_{1}+\cdots+l_{k}=j}\binom{j}{l_{1},\dots,l_{k}}r^{m-j}E\bigg[S_{k}^{n-i}\prod_{i=1}^{k}Y_{i}^{l_{i}}\bigg],\nonumber
\end{align}
where $l_{1},l_{2},\dots,l_{k}$ are positive integers.\\
Therefore, by \eqref{25}, we obtain the following theorem.
\begin{theorem}
For $m,n\ge 0$, we have
\begin{align*}
&\phi_{m+n,r}^{Y}(x)\\
&= \sum_{i=0}^{n}\sum_{k=0}^{m}\binom{n}{i}\phi_{i,r}^{Y}(x)x^{k}\frac{1}{k!}\sum_{j=k}^{m}\binom{m}{j} \sum_{l_{1}+\cdots+l_{k}=j}\binom{j}{l_{1},\dots,l_{k}}r^{m-j}E\bigg[S_{k}^{n-i}\prod_{i=1}^{k}Y_{i}^{l_{i}}\bigg],
\end{align*}
where $l_{1},l_{2},\dots,l_{k}$ are positive integers.
\end{theorem}
From Theorem 2.5, we get
\begin{align}
\phi_{m+n,r}^{Y}(x)&=\lim_{y\rightarrow 0}\phi_{m+n,r}^{Y}\bigg(\frac{x}{y},y\bigg) \label{26}\\
&=\lim_{y\rightarrow 0}\sum_{i=0}^{n}\sum_{k=0}^{m}\binom{n}{i}\bigg(\frac{x}{y}\bigg)_{k}y^{k}\phi_{i}^{Y}\bigg(\frac{x}{y}-k,y\bigg)\frac{1}{k!} \sum_{j=k}^{m}\binom{m}{j} \nonumber\\
&\quad\quad\quad\quad \times \sum_{l_{1}+\cdots+l_{k}=j}\binom{j}{l_{1},\dots,l_{k}}E\bigg[(S_{k}+r)^{n-i}\prod_{i=1}^{k}Y_{i}^{l_{i}}\bigg]r^{m-j}\nonumber \\
&=\sum_{i=0}^{n}\sum_{k=0}^{m}\binom{n}{i}x^{k}\phi_{i}^{Y}(x)\frac{1}{k!}\sum_{j=k}^{m}\binom{m}{j} \sum_{l_{1}+\cdots+l_{k}=j}\binom{j}{l_{1},\dots,l_{k}}E\bigg[(S_{k}+r)^{n-i}\prod_{i=1}^{k}Y_{i}^{l_{i}}\bigg]r^{m-j},\nonumber
\end{align}
where $l_{1},l_{2},\dots,l_{k}$ are positive integers. \\
Therefore, by \eqref{26}, we obtain the following theorem.
\begin{theorem}
For $m,n\ge 0$, we have
\begin{displaymath}
\phi_{m+n,r}^{Y}(x)=\sum_{i=0}^{n}\sum_{k=0}^{m}\binom{n}{i}x^{k}\phi_{i}^{Y}(x)\frac{1}{k!}\sum_{j=k}^{m}\binom{m}{j}\sum_{l_{1}+\cdots+l_{k}=j}\binom{j}{l_{1},\dots,l_{k}} E\bigg[(S_{k}+r)^{n-i}\prod_{i=1}^{k}Y_{i}^{l_{i}}\bigg]r^{m-j},
\end{displaymath}
where $l_{1},l_{2},\dots,l_{k}$ are positive integers.
\end{theorem}
When $Y=1$, the expression in Theorem 2.7 boils down to
\begin{align*}
\phi_{m+n}(x)&= \sum_{i=0}^{n}\sum_{k=0}^{m}\binom{n}{i}x^{k}\phi_{i}(x)\frac{1}{k!}\sum_{j=k}^{m}\binom{m}{j}\sum_{l_{1}+\cdots+l_{k}=j}r^{m-j}(k+r)^{n-i} \\
&=\sum_{i=0}^{n}\sum_{k=0}^{m}\binom{n}{i}x^{k}\phi_{i}(x){m+r \brace k+r}_{r}(k+r)^{n-i}.
\end{align*}
\section{Conclusion}
Let $Y$ be a random variable such that the moment generating function of $Y$ exists in a neighborhood of the origin. We introduced probabilistic extensions of the bivariate Bell polynomials and the bivariate $r$-Bell polynomials, namely the probabilistic bivariate Bell polynomials associated with $Y$ and the probabilistic bivariate $r$-Bell polynomials associated with $Y$. Zheng and Li found recurrence relations for the bivariate Bell polynomials and the bivariate $r$-Bell polynomials. In this paper, by utilizing the simple fact \eqref{11} we obtained recurrence relations for the probabilistic bivariate Bell polynomials associated with $Y$ and the probabilistic bivariate $r$-Bell polynomials associated with $Y$. These relations boil down to the
relations found by Zheng and Li when $Y=1$.\par
As one of our future projects, we would like to continue to study probabilistic extensions of many special polynomials and numbers and to find their applications to physics, science and engineering as well as to mathematics.

\section{Declarations}
The authors declare that they have no conflicts of interest to report regarding
the present study.

\end{document}